\documentclass[amsart]{article}

\usepackage{setspace}
\usepackage{amsmath}
\usepackage{amssymb}
\usepackage{amscd}
\usepackage{bbm}
\usepackage{cyrillic}
\usepackage{dirtytalk}
\usepackage{latexsym}
\usepackage{makeidx}
\usepackage{mathrsfs}
\usepackage{mathtools}
\usepackage{theorem}
\usepackage{color}
\usepackage[all]{xy}
\usepackage[colorlinks]{hyperref}
\hypersetup{colorlinks=false , linkcolor=blue}

\DeclareMathAlphabet{\mathbbold}{U}{bbold}{m}{n}

\newtheorem{teo}{Theorem}[section]

\newtheorem{lema}[teo]{Lemma}
\newtheorem{coro}[teo]{Corollary}

{\theorembodyfont{\rmfamily }
\newtheorem{defi}[teo]{Definition}

\newtheorem{obser}[teo]{Remark}
\newtheorem{ejem}[teo]{Example}
\newtheorem{parra}[teo]{}}
\def\demo{\noindent \textit{Proof: }}
\def\fin{\end{document}}

 \DeclareMathOperator{\hofib}{hofib}
\DeclareMathOperator{\Hom}{Hom}

 \DeclareMathOperator{\Spec}{Spec}
\DeclareMathOperator{\Td}{Td} \DeclareMathOperator{\Th}{Th}

\newcommand{\cyrrm}{\fontencoding{OT2}\selectfont\textcyrup}
\def\11{1}
\def\11{\mathbbm{1}}

\def\ad{\mathrm{ad}}

\def\BM{\mathrm{BM}}
\def\cf{\emph{cf. }}
\def\ch{\mathrm{ch}}

\def\hatch{\widehat{\mathrm{ch}}}

\def\hatKH{\widehat{KH}}

\def\p{\mathfrak{p}}

\def\raya{\ \underline{\phantom{a}}\ }
\def\refi{\mathrm{ref}}

\def\t{\mathbf{t}}

\def\qed{\hspace*{\fill }$\square $ }

                    \def\AA{\mathbb{A}}

 \def\EE{\mathbb{E}}
 \def\FF{\mathbb{F}}

\def\M{\mathcal{M}} \def\MM{\mathbb{M}}
\def\O{\mathcal{O}}
 \def\PP{\mathbb{P}}
                    \def\QQ{\mathbb{Q}}

 \def\ZZ{\mathbb{Z}}

\def\KGL{\mathrm{KGL}}

\def\HB{\mathrm{H}_{\cyrrm{B}}}
\def\hatHB{\widehat{\mathrm{H}}_{\cyrrm{B}}}
\def\hatKGL{\widehat{\mathrm{KGL}}}

\def\SH{\mathbf{SH}}

\begin{document}

\title{On the Riemann-Roch formula  without projective hypothesis}

\author{A. Navarro \thanks{Institut f\"{u}r Mathematik (U. Z\"{u}rich)},
J. Navarro \thanks{Departamento de Matem\'{a}ticas, Universidad de Extremadura, Spain} }


\maketitle

%
%

\begin{abstract}

Let $S$ be a finite dimensional noetherian scheme. For any proper morphism between
smooth $S$-schemes, we prove a Riemann-Roch formula relating higher algebraic
$K$-theory and motivic cohomology, thus with no projective hypothesis neither on the
schemes nor on the morphism. We also prove, without projective
assumptions, an arithmetic Riemann-Roch theorem involving Arakelov's higher $K$-theory and
motivic cohomology as well as an analogue result for the relative cohomology of a
morphism.

These results are obtained as corollaries of a motivic statement that is valid for morphisms
between oriented absolute spectra in the stable homotopy category of $S$.

\end{abstract}



\section*{Introduction}

Let $\,f\colon Y\to X\,$  be a proper morphism between nonsingular quasi-projective
varieties over a field $k$. The original Grothendieck's Riemann-Roch theorem
(\emph{cf.} \cite{BorelSerre}) states that the square
$$
\xymatrix{
K_0(Y)\ar[r]^{f_!} \ar[d]_{\Td(T_Y)\ch}& K_0(X)\ar[d]^{\Td(T_X)\ch}\\
CH^\bullet (Y) _\QQ \ar[r]^{ f_*} & CH^\bullet (X) _\QQ }
$$
commutes, or, in other words, that for any $a \in K_0 (Y)$ the formula
\begin{equation}\label{RRformula}
\Td (T_{\scriptscriptstyle X})\, \ch (f_!(a)) \, =\,  f_* \bigl( \Td
(T_{\scriptscriptstyle Y})\, \ch (a)\bigr) \
\end{equation}
holds. This diagram relates the exceptional direct image $f_!$ on the Grothendieck
group of vector bundles $K_0$ to the direct image $f_*$ on the Chow ring with
rational coefficients. The vertical arrows involve the Chern character $\ch$,
adequately twisted with the Todd series $\Td$ of the corresponding tangent bundles.

Later on, Grothendieck aimed to generalise the Riemann-Roch theorem in three
directions: allowing general schemes not necessarily over a base field, replacing the
smoothness  condition on the schemes by a regularity condition on the morphism, and
avoiding any projective assumption either on the morphism or on the schemes. The
first two directions were one of the main objectives of \cite{SGA6} (\cf
\cite[O.1]{SGA6}). However, the third one---namely, proving the Riemann-Roch for
proper complete intersection morphisms between (non-projective) schemes---remained as
an open question proposed by Grothendieck at \cite[XIV.2]{SGA6} as \emph{the
Riemann-Roch without projective hypothesis}.



In subsequent years there were results in this direction. If $X$ is a smooth and
proper variety over a field $k$ of characteristic zero, the Chow lemma and the
resolution of singularities allow to construct a projective variety $\overline{X}$
and a birrational map $\overline{X} \to X$ that is a composition of blowing ups along
smooth centers. Fulton used this idea in \cite{Fulton2} to prove the Riemann-Roch
formula for the projection $X\to \Spec (k)$. In fact, the general case of a proper
morphism between smooth varieties over a field of characteristic zero also follows
using this idea (\cf \cite{Papa}, \cite{FG}).

The case of positive characteristic was proved in \cite{FG} using a different
approach. Fulton-Gillet's  proof avoids the resolution of singularities due to the
formulation of the Riemann-Roch theorem developed in \cite{BFM}, the technique of
envelopes (\cite{Gillet2}) and a little use of higher $K$-theory.

The result of \cite[Corollary 2]{FG} we just mentioned is the furthest generalization
we know of the Riemann-Roch without projective hypothesis. Nevertheless, despite the
Riemann-Roch for higher $K$-theory was known at that moment (\cf \cite{Gillet}),
Fulton-Gillet's proof did not extend to higher $K$-theory and they explicitly asked
about that case (\cf \cite[3.1.4]{FG}). Note that their proof also assumed schemes to
be over a field. As Grothendieck already pointed out,
\begin{quotation}
\say{Il semble clair que la d\'{e}monstration de la formule Riemann-Roch dans ce cas
g\'{e}n\'{e}ral [i.e., without projective hypothesis] demandera l'introduction d'id\'{e}es
essentiellement nouvelles}(\cite[XIV.2]{SGA6})
\end{quotation}
and every proof restricts to the case where its new ideas may be applied.








In this paper we prove a higher Riemann-Roch theorem without projective hypothesis,
neither on the schemes, nor on the morphism. Our statement is valid for smooth
schemes defined over a finite dimensional, noetherian scheme $S$. To be more precise,
for any smooth $S$-scheme $\,X$, let us write $\ch \colon K(X) \to H_{\M} (X, \QQ)$
for the higher Chern character, relating higher $K$-groups to motivic cohomology
groups with rational coefficients (see definitions in Example \ref{Ejem espectros}):

\medskip

\goodbreak

\noindent \textbf{Theorem \ref{Teo Riemann Roch}} {\it Let $S$ be a finite
dimensional noetherian scheme, and let $f\colon Y\to  X$ be a proper morphism between
smooth $S$-schemes. The following diagram commutes:
$$
\xymatrix{
K(Y)\ar[r]^{f_*} \ar[d]_{\Td(T_Y)\ch}& K(X)\ar[d]^{\Td(T_X)\ch}\\
H_{\M}(Y,\QQ) \ar[r]^{ f_*} & H_{\M}(X,\QQ) }
$$
In other words, for any element $\,a \in K(Y)$, the following formula holds:
$$
\Td (T_{\scriptscriptstyle X})\, \ch (f_*(a)) \, =\,  f_* \bigl( \Td
(T_{\scriptscriptstyle Y})\, \ch (a)\bigr) \ .
$$
}

\medskip

Note that the restriction of the formula above to the $K_0$ group produces a
Grothendieck-Riemann-Roch statement that is valid for proper morphisms between
general (non-projective) smooth $S$-schemes.

In addition, we also prove other Riemann-Roch statements for the relative cohomology
of a morphism (\cf Theorem \ref{Teo RR Relativo}) and for the arithmetic counterparts
of higher $K$-theory and motivic cohomology (\cf Theorem \ref{Teo RR Aritmetico}).
These results also improve previous versions of \cite{Navarro Articulo} and
\cite{HS}, respectively.



Our proof relies on a duality theorem of Ayoub (\cite{Ayoub}), as refined by
Cisinski-D\'{e}glise in \cite{CD}, and on the Riemann-Roch formula for closed embeddings.
Roughly speaking, since duality holds without projective assumptions, the
Riemann-Roch formula can also be proved in this generality. The main idea---namely,
that the Riemann-Roch theorem is a comparison of dualities---was already pointed out
by D\'{e}glise at his lecture at the $K$-theory satellite conference of ICM 2014 at
Beijing (\cite{Deglise4}).

The paper is written in Morel-Voevodsky's language of $\AA^1$-homotopy theory
(\cite{MV}, \cite{Voevodsky}). In addition, our development of the Riemann-Roch is
directly inspired by Panin's work on orientation theory (\cite{Panin5},
\cite{Panin4}, \cite{Panin2}, \cite{Navarro Grothendieck}), that has been lifted to
the motivic homotopy setting by D\'{e}glise (\cite{Deglise}).

Finally, let us briefly comment the plan of this article. The first, preliminary
section, recalls definitions and statements from stable homotopy, in particular
regarding cohomology, orientations, Borel-Moore homology and Thom spaces. In Section
2 we construct a funtorial direct image in cohomology as the dual of the direct image
in Borel-Moore homology, and prove its main properties. Finally, in the last section
we establish the aforementioned Riemann-Roch theorems.

\subsection*{Acknowledgements}

The authors would like to thank J. Ayoub, J.I. Burgos Gil and F. D\'{e}glise for many
helpful discussions and comments. The first author has been partially supported by
ICMAT Severo Ochoa project SEV-2015-0554 (MINECO), MTM2013-42135-P (MINECO) and
ANR-12-BS01-0002 (ANR). The second author has been partially supported by Junta de
Extremadura and FEDER funds.

\section{Preliminaries}

All schemes we will consider throughout this paper are smooth over a finite
dimensional noetherian base $S$. We use the same notations of \cite{Navarro Articulo}
and recall the indispensable.

\medskip

Let $X$ be a smooth $S$-scheme and denote $\SH(X)$ the stable homotopy category of
Voevodsky, whose objects are called \emph{spectra} (\cf \cite{Voevodsky}).

These categories are symmetric monoidal, meaning that there is a symmetric product,
$\wedge$, and a unit object $\11_X$. There is also a shift functor $ [1]\colon\SH
(X)\to  \SH (X) $ and a Tate object $\11_X(1)$, that defines a twist $ \raya \wedge
\11_X(1)\colon\SH (X)\to  \SH (X) $. If $E$ is an spectrum of $\SH(X)$ and $p,q\in
\ZZ$, we write $E(q)[p]$ for the result of shifting $p$ times and twisting $q$ times.

Moreover, the family of categories obtained if one varies $X$ satisfies
Grothendieck's six functors formalism (\cf \cite{Ayoub} and \cite{CD}). In
particular, for any morphism $f\colon Y\to X$ there exists a pair of adjoint functors
$$
f^* \colon \SH (X) \leftrightarrows \SH (Y) \colon f_*.
$$

If $f$ is separated of finite type, then there also exist mutually adjoint
exceptional functors
$$
f_! \colon \SH (Y) \leftrightarrows \SH (X) \colon f^! \ .
$$

Finally, if $\pi \colon X\to T$ is a smooth morphism, then there
exists an ``extension by zero" functor $\pi_\sharp \colon \SH (X) \to \SH (T)$ that
is left adjoint to $\pi^*$. These functors satisfy adequate functorial properties
(\cf \cite[1.4.2]{Ayoub} and the refinement in \cite[A.5]{CD}).

\begin{defi}
A \textit{ring spectrum} is an associative commutative unitary monoid in the stable
homotopy category $\SH (X)$ of a scheme $X$.

An \textit{absolute ring spectrum} $\EE$ is a ring spectrum on the stable homotopy
category of the base scheme $\SH(S)$. Equivalently, an absolute ring spectrum is a
family of ring spectra $\EE_X$ on $\SH(X)$ for every scheme $X$ that is stable by
pullback---i.e., such that for every morphism of schemes $f\colon Y\to X$ we have
fixed an isomorphism $\epsilon_f\colon f^*\EE_X \to \EE_Y$ satisfying the usual
compatibility conditions (\cf \cite[1.2.1]{Deglise}).

A \textit{morphism} of absolute ring spectra $\varphi \colon \EE \to \FF$ is a
morphism of spectra in $\SH(S)$ or, equivalently, a collection of morphisms of ring
spectra $\varphi_X \colon \EE_X \to \FF_X$, one for each scheme $X$, that is stable
by pullback.

If $\EE$ is an absolute ring spectrum, the {\it $\EE$-cohomology} of a scheme $X$ is,
for $p,q \in \mathbb{Z}$:
$$
\EE^{p,q} (X) :=  \Hom_{\SH(X)} (\11_X ,\EE_{\scriptscriptstyle X}(q)[p]) \quad
 \quad , \hspace{.5cm} \EE (X) \coloneqq \bigoplus_{p,q}
\EE^{p,q}(X) \ .
$$

An {\it oriented absolute ring spectrum} is a pair $(\EE , c_1)$, where $\EE$ is an
absolute ring spectrum and $c_1 \in \EE^{2,1} (\mathbb{P}^\infty)$ is a class
satisfying that $i_1^* (c_1) = \eta$, for $i_1 \colon \mathbb{P}_1 \hookrightarrow
\mathbb{P}^\infty$ and $\eta \in \EE^{2,1} (\mathbb{P}_1)$  the canonical class
defined by the Tate object (\cf \cite[\S 1.3]{Navarro Articulo}).

\end{defi}

\begin{ejem}\label{Ejem espectros}
There is a pleiad of cohomology theories represented by oriented absolute ring
spectra. Let us recall here those that we will mention later on, although the reader
may consult \cite[1.4]{Navarro Articulo} for more examples to whom apply the results
of this paper.

\begin{itemize}
\item
The {\it $K$-theory spectrum} $\KGL$ is an absolute ring spectrum defined in
\cite{Voevodsky} and \cite{Riou3}.
It represents Weibel's homotopy invariant $K$-theory (\cf
\cite[2.15]{Cisinski})---and therefore also represents Quillen's algebraic $K$-theory
for regular schemes. The cohomology it defines will be denoted $K(\raya )$.

\item
The $\QQ$-localization of the $K$-theory spectrum admits
a decomposition induced by the Adams operations, i.e., $ \KGL_\QQ = \bigoplus_{i\in
\ZZ} \KGL_ \QQ ^{(i)}\in \SH(S)_\QQ, $ where $\KGL_ \QQ ^{(i)}$ denotes the
eigenspaces for the Adams operations (\cite[5.3]{Riou3}).

Beilinson's {\it motivic cohomology spectrum} is defined as $\HB=\KGL_\QQ^{(0)} \in
\SH(S)_\QQ$ and it is also an absolute ring spectrum. It represents motivic
cohomology with rational coefficients, which we denote
$H_{\mathcal{M}}^{p}(\raya,\QQ(q))$.

\item
Let $g\colon T \to X$ be a morphism of schemes and let $\EE$ be an absolute ring
spectrum.

The spectrum $\hofib_\EE (g)\coloneqq \hofib (\EE_X\to g_*g^*\EE_X)$ defines
cohomology groups that are called the \emph{relative cohomology} of $g$:
$$\EE^{*,*}(g)\coloneqq
\Hom_{\SH(X)}(\11_X,\hofib(g)[*](*)) \ . $$

They fit into long exact sequences
$$
\cdots \to \EE^{p,q}(g)\to \EE^{p,q}(X)\to \EE^{p,q}(T)\to \EE^{p+1,q}(g)\to\cdots.
$$
and generalize many constructions: the cohomology with proper support, the cohomology
with support on a closed subscheme and the reduced cohomology are, respectively, the
relative cohomology of a closed immersion, an open immersion and the projection over
a base point (\cf \cite[1.19]{Navarro Articulo}).

Let $f\colon Y\to X$ be another morphism and denote $g_{\scriptscriptstyle Y}\colon
T\times_X Y\to Y$. If either $g$ is proper or $f$ is smooth, then the absolute
spectrum over $X$ defined by $\hofib_\EE(g)$ represents $\EE(g_{\scriptscriptstyle
Y})$ (\cf \cite[\S 1.17]{Navarro Articulo}).
\end{itemize}
\end{ejem}

\begin{defi}
Let $E$ in $\SH(X)$ be a ring spectrum. An \textit{$E$-module} is a spectrum $M$ in
$\SH(X)$ together with a morphism of spectra $\upsilon\colon E\wedge M\to M$ in
$\SH(X)$ satisfying the usual module condition. An \textit{absolute $\EE$-module} is
a module over an absolute ring spectrum $\EE$.

Let $\varphi\colon E\to F$ be a
morphism of ring spectra and $(M' , \, v')$ be an $F$-module. A
\textit{$\varphi$-morphism} of modules $\Phi \colon M\to M'$ is a morphism of spectra
in $\SH(X)$ such that $ v'\circ( \varphi \wedge \Phi)=\Phi\circ v$.
\end{defi}

\begin{ejem}\label{Ejem modulos}
Holmstrom-Scholbach defined in \cite{HS} the {\it arithmetic $K$-theory} and the {\it
arithmetic motivic cohomology} ring spectra, $\hatKGL$ and $\hatHB$. They are
absolute modules over $\KGL$ and $\HB$, respectively. In addition, there is an
arithmetic Chern character $\hatch\colon \hatKGL\to \hatHB$, that is a $\ch$-morphism
of modules (\cf \cite[4.2]{HS}).
\end{ejem}

\begin{obser}\label{Obser olvidar soporte}
The machinery of the six functors formalism provides the $\EE$-cohomology defined by
oriented absolute spectra with the classic properties of cohomologies (such as
inverse image, cup product, Chern classes, etc.). For a complete account of these
properties the reader may consult \cite{Deglise}.

Let us recall here the construction of the
morphism of \emph{forgetting support}, since it will be used later on. Let
$Z\xrightarrow j Y \xrightarrow i X$ be closed immersions and let $\MM$ be an absolute
module. We define a morphism
$$
j_\flat  \colon\MM_Z (X) \to \MM_Y(X)
$$
as follows: the adjunction of $(j^* ,\, j_*)$ induces a morphism $i_*(\ad) \colon i_*
(\11_Y)\to i_*(j_*j^*\11_Y)=(ij)_*\11_Z$. Let $a\colon (ij)_*\11_Z\to \MM_X $ be in
$\MM_Z (X)$. The element $j_\flat(a)\in \MM_Y(X)$ is defined as:
$$
i_* (\11_Y)\xrightarrow {i_*(\ad)}(ij)_*\11_Z\xrightarrow
{a}\MM_X \ .
$$
\end{obser}

\medskip
\begin{parra}[Borel-Moore homology]
Let $\EE$ be an absolute ring spectrum, and $\MM$ be an $\EE$-module. The
\textit{Borel-Moore homology} of a scheme $\pi_{\scriptscriptstyle X} \colon
X\rightarrow S$ is defined as
$$
\MM^\BM_{p,q}(X) \coloneqq \Hom _{\SH(X)}(\11_{\scriptscriptstyle X},
\pi_{\scriptscriptstyle X}^!\MM_{\scriptscriptstyle S}(-q)[-p]) \hspace{.3cm} ,
\hspace{.3cm} \MM^\BM(X) \coloneqq \bigoplus_{p,q} \MM^\BM_{p,q}(X) \ .
$$
Observe the equivalent descriptions, using adjunction:
\begin{align*}\label{HomologyAdj}
\MM^\BM_{p,q}(X) &=   \Hom_{\SH(S)} (\pi_{\scriptscriptstyle
X!} \11_{\scriptscriptstyle X}, \MM_{\scriptscriptstyle S}(-q)[-p]))  \\
&= \Hom _{\SH(S)}(\11_{\scriptscriptstyle S},\pi_{\scriptscriptstyle
X*} \pi_{\scriptscriptstyle X}^!\MM_{\scriptscriptstyle S}(-q)[-p]) \ .
\end{align*}

\medskip

If $f\colon Y\to X$ is a proper morphism, then the direct image with compact support,
$f_!$, is canonically isomorphic to $f_*$ (\cite[2.2.8]{CD}), so that the adjunction
of $(f_!,f^!)$ induces a natural transformation
$$
\pi_{{\scriptscriptstyle Y}*} \pi_{{\scriptscriptstyle Y}}^! \, =\,
\pi_{{\scriptscriptstyle X}*} f_* f^! \pi_{\scriptscriptstyle X}^! \, = \,
\pi_{{\scriptscriptstyle X}*} f_!f^!\pi_{\scriptscriptstyle X}^! \, \buildrel{\,
\mathrm{ad} \, }\over{\longrightarrow} \, \pi_{{\scriptscriptstyle
X}*}\pi_{\scriptscriptstyle X}^! \ ,
$$
which in turn produces a {\it direct image} in BM-homology
$$
f_*\colon \MM^\BM_{p,q}(Y)\longrightarrow \MM^\BM_{p,q}(X).
$$

As expected Borel-Moore homology is a module, in the classic sense, over the
$\EE$-cohomology. More precisely, there is a \emph{cap product}:
$$
\EE^{p,q}(X)\times\EE^\BM_{r,s}(X) \longrightarrow  \EE^{\BM}_{r-p, s-q}(X),
$$
where a pair $(a,m)$ is mapped into the element $a\cdot m \in \EE^{\BM}_{r-p,
s-q}(X)$ defined as:
\begin{multline*}
\11_X\xrightarrow{a} \EE_X(q)[p]\wedge\11_X\xrightarrow{m}
\pi_X^*(\EE_S)\wedge \pi_X^!(\EE_S)(q-s)[p-r]\xrightarrow{\sim}\\
\xrightarrow{\sim}\pi_X^!(\EE_S\wedge\EE_S)(q-s)[p-r]\xrightarrow{\mu}\pi_X^!
(\EE_S)(q-s)[p-r].
\end{multline*}
Here, $\mu$ stands for the structural map of
$\EE$ and the first isomorphism is due to the fact that the structural map of the
bilateral module \footnote{For the definition of left and right module in this
context see \cite[2.1.94]{Ayoub}.} $[\pi^*, \pi^!]$
is an isomorphism (\cf \cite[2.3.27]{Ayoub}).
\end{parra}

\medskip

\begin{parra}[The Thom space]
Let $V\to X$ be a vector bundle of rank $d$ and let us write $\bar V \coloneqq
\PP(V\oplus \11)$ for the projective completion.

The \textit{Thom space} of $V$ is defined as:
$$
\Th (V)\coloneqq V/V-\{0\}\simeq \bar V/\PP(V) \in \SH(X) \ .
$$

Let $(\EE,c_1)$ be an oriented absolute ring spectrum. Let us write
$$
\EE^{*,*} (\Th (V) ) \coloneqq \Hom_{\SH(X)} (\Th(V),\EE_X(*)[*]).
$$
There exists a long exact sequence
$$
\cdots \to \EE^{*,*}(\Th(V))\xrightarrow{\pi}\EE^{*,*}(\bar V)\xrightarrow{i^*}
\EE^{*,*}(\PP(V))\to \cdots.
$$

The \textit{Thom class} is the following element in cohomology
$$
\t (V) \coloneqq \sum _{i=0}^d (-1)^i c_i (V)x^i \  \in \, \EE^{2d,d}(\bar V) \ \ \ ,
\hspace{.5cm} x\coloneqq c_1(\O_{\bar V}(-1))\, .
$$
As a consequence of the projective bundle theorem (\cite[2.1.13]{Deglise}), $i^*(\t
(V))=0$, so that there exists a unique class $\t^{\refi} (V)\in \EE^{2d,d}(\Th(V))$,
called the \textit{refined Thom class}, such that $\pi(\t^\refi(V))=\t(V)$. Moreover,
the cohomology of the Thom space $\EE(\Th(V))$ is a free $\EE(X)$-module of rank one
generated by $\t ^\refi (V)\in \EE^{2d,d}(\Th(V))$.

Also, let us recall that the deformation to the normal bundle induces an
isomorphism
\begin{equation}\label{Eq iso Thom cohomologia con soporte}
\begin{matrix}
\EE(\Th(V))&\buildrel{\sim}\over{\longrightarrow }& \EE_X(\bar V)\\
\t ^{\refi}(V)& \mapsto & \bar \eta_X^V,
\end{matrix}
\end{equation}
where $\bar \eta_X^V$ denotes the refined fundamental class of $X$ on $V$ (\cf
\cite[2.3.1]{Deglise}).

\medskip
Similar statements are true for $\EE$-modules. Due to the lack of a reference, let us
sketch the proof of the projective bundle theorem in this case:

\begin{lema}\label{Lema Fibrado Proyectivo Modulos}
Let $(\EE,c_1)$ be an absolute ring spectrum, $\MM$ be an absolute $\EE$-module, and
$V\to X$ be a vector bundle.

It holds:
$$
\MM(\PP(V))= \MM(X) \otimes_{\EE(X)} \EE(\PP(V)) \ .
$$
\end{lema}

\demo As in \cite[3.2]{Deglise2}, due to Mayer-Vietoris we can reduce to the case of
a trivial vector bundle. To prove $\MM(\PP^n_X)= \MM(X) \otimes_{\EE(X)}
\EE(\PP^n_X)$ we can assume $X=S$ and proceed by induction on $n$. The case $n=0$ is
trivial; for the induction step consider the homotopy cofiber sequence
$$
\PP^{n-1}\xrightarrow i \PP^n\xrightarrow \pi \PP^n/\PP^{n-1}
$$
that produces a long exact sequence
$$
\cdots \to \MM(\PP^n/\PP^{n-1})\to \MM(\PP^n)\xrightarrow
{i^*}\MM(\PP^{n-1})\to\cdots.
$$
Since $\PP^n/\PP^{n-1}\simeq \11(n)(2n)$ (\cf \cite[3.2.18]{MV}) then
$\MM^{p,q}(\PP^n/\PP^{n-1})=\MM^{p-2n,q-n}(S)$. By induction hypothesis we deduce
$i^*$ is surjective (since we have $i^*(m\cdot c_1(\O_{\PP^n}(-1))^j)=m\cdot
c_1(\O_{\PP^{n-1}}(-1))^j)$ for $m\in \MM(X)$). The map $\pi\colon
\MM^{p-2n,q-n}(S)\to \MM^{p,q}(\PP^n)$ maps $m\mapsto m\cdot
c_1(\O_{\PP^{n-1}}(-1))^n$ (\cf \cite[Proof 2.1.3]{Deglise}). Hence, the thesis
follows.

\qed

\medskip
As a consequence of this Lemma, for any vector bundle $\,V \to X$, the
$\MM$-cohomology of its Thom space is isomorphic to the $\MM$-cohomology of $\,X$
through the morphism
\begin{equation}\label{Ec M-cohomologia Thom}
\begin{matrix}
\MM(X) &\buildrel{\sim}\over{\longrightarrow }&  \MM^{p,q}(\Th(V))\\
m& \mapsto & m\cdot \t ^{\refi}(V).
\end{matrix}
\end{equation}

\end{parra}

\begin{obser}
More generally, all these statements above remain true if $V$ is a virtual vector
bundle. Indeed, as shown by Riou (\cite{Riou3}), the Thom space construction in $\SH$
extends to a canonical functor
$$
\Th\colon \underline{K(X)} \to \SH (X)
$$
defined on the category $\underline{K(X)}$ of virtual vector bundles
(\cite[4.12]{Del87}).

In particular, this means that for every short exact sequence of vector bundles $0\to
E'\to E \to E''\to 0$ in $X$ there is a canonical isomorphism in $\SH (X)$:
\begin{equation}\label{Thom y sucesion exacta}
\Th(E')\wedge \Th(E'')\xrightarrow \sim \Th (E) \ .
\end{equation}
It induces an isomorphism (\cf \cite[2.4.8]{Deglise}):
$$
\begin{matrix}
\EE (\Th(E'))\otimes \EE (\Th(E'')) & \buildrel{\sim}\over{\longrightarrow} &
\EE(\Th (E)),\\
\t^\refi (E')\otimes \t^\refi(E'') & \mapsto & \t^\refi(E')\cdot \t^\refi(E'')=
\t^\refi (E)
\end{matrix}
$$

Thus, the Thom class of a virtual vector bundle
$\xi=[E]-[E']$ is defined to be the class
$$
\t^\refi(\xi)=\t^\refi(E)\cdot \t^\refi(E')^{-1} \ .
$$
\end{obser}

\section{Direct image for proper morphisms}

In what follows, we will make a strong use of the following statement, which is a
particular case of a general duality theorem of Ayoub (\cf \cite[1.4.2]{Ayoub},
\cite[2.4.28]{CD}).

\begin{teo}\label{Ayoub}
If $\pi \colon X\to T$ is a smooth morphism, there exists a functorial isomorphism
$$
\pi^* (\raya) \wedge \Th(T_{\pi}) \buildrel{\ \ \sim \ \ }\over{\longrightarrow }
\pi^! (\raya )  ,
$$
where $T_{\pi}$ stands for the virtual tangent bundle of $\pi$ (\cf
\cite[B.2.7]{Fulton}).
\end{teo}
\qed

\bigskip

In the particular case of a smooth $S$-scheme, $\pi_{\scriptscriptstyle X} \colon
X\to S$, the tangent bundle of $\pi_{\scriptscriptstyle X}$ is the tangent bundle
$T_{\scriptscriptstyle X}$ of $X$. Hence, the theorem above implies the following
isomorphism:
$$
\Hom_{\SH(X)} (\Th(-T_{\scriptscriptstyle X}) ,\MM_{\scriptscriptstyle
X}(q)[p])\simeq \Hom_{\SH(X)} (\11_{\scriptscriptstyle X},\pi_{\scriptscriptstyle
X}^! \MM_{\scriptscriptstyle S}(q)[p]))
$$
that relate the cohomology of the Thom space of the tangent bundle to its Borel-Moore
homology; that is to say,
\begin{equation}\label{ThomYHomologia}
\MM^{p,q}(\Th (-T_{\scriptscriptstyle X}))\ \simeq \ \MM^{\BM}_{-p,-q}(X) \ .
\end{equation}

Combining these facts with the computation of the cohomology of the Thom space
(\ref{Ec M-cohomologia Thom}), we obtain duality isomorphisms between the
$\MM$-cohomology and the Borel-Moore homology of a smooth scheme:
\begin{equation} \label{IsosDualidad}
\MM^{p,q} (X) \ \simeq \ \MM^{p+2n,q+n} ( \Th (-T_X) )
\ \simeq \ \MM^{\BM}_{-p-2n , -q -n} (X) \ .
\end{equation}

With these isomorphisms, it is possible to define a direct image in cohomology
without projective assumptions:

\begin{defi}\label{LaDefinicion}
Let $(\EE, c_1)$ be an absolute oriented ring spectrum and $\MM$ be an absolute
$\EE$-module. The \textit{direct image}---in $\MM$-cohomology---of a proper morphism
$f\colon Y \to X$ between smooth $S$-schemes is defined as the composition
\begin{equation*}
\MM^{p,q}(Y) \, \simeq \, \MM^{\BM}_{-p-2n,-q-n}(Y) \, \buildrel{\, f_* \,
}\over{\longrightarrow} \, \MM^{\BM}_{-p-2n,-q-n}(X) \, \simeq \, \MM^{p+2d,q+d}(X)\
,
\end{equation*}
where $n=\dim_S Y$ and $d=n-\dim_SX$.
\end{defi}

\begin{obser}\label{Obser Gysin Thom}
Due to Ayoub's duality, it is not surprising that we can express the direct image in
terms of Thom spaces. Let us describe in detail the case of a closed immersion
$i\colon Z \to X$.

Let $N_{Z/X}$ be the normal bundle and recall $T_Z=i^*T_X-N_{Z/X}\in K_0(Z)$. Then
$$
\Th(T_Z)=\Th(-N_{Z/X})\wedge \Th(i^*T_X)=i^!(\11_X)\wedge i^*(\Th(T_X)\, ,
$$
because of isomorphism (\ref{Thom y sucesion exacta}) and the fact that
$i^!(\11_X)=\Th(-N_{Z/X})$ (\cf \cite[1.3.5.3]{Deglise}) and
$\Th(i^*T_X)=i^*(\Th(T_X))$.

Therefore, there is a map
$$
\xymatrix@R=6pt{\EE (\Th(-T_Z)) \ar@{=}[d] \ar@{-->}[rr] & & \EE(\Th(-T_X))
\ar@{=}[dd] \\
\Hom_{\SH(Z)}(\11_Z, \Th(T_Z)\wedge \EE_Z) \ar@{=}[d]& & \\
\Hom_{\SH(X)}(\11_X, i_!i^!(\11_X)\wedge \Th (T_X)\wedge
\EE_X)\ar[rr]^-{\mathrm{ad}_{i_!i^!}} & & \Hom_{\SH(X)} (\Th(-T_X), \EE_X), }
$$
where $\mathrm{ad}_{i_!i^!}$ is obtained by composing with the natural adjunction
$i_!i^!\to \mathrm{Id}$.

Denote $\pi=\pi_X$ and $\pi'=\pi_Z$. Observe that $T_{\pi i }=T_{\pi' }=T_Z$, then
$\Th(T_{\pi i})=\Th(T_Z)= i^!(\11_X)\wedge i^*(\Th(T_{\pi})$. The compatibility of
the above map $\EE (\Th(-T_Z))\to \EE (\Th(-T_X))$ and the direct image in
Borel-Moore homology derives from the commutative square
$$
\xymatrix@R=6pt{\pi_! i_! i^! \pi^! (\raya ) \ar@{=}[d]
\ar[rr]^-{\mathrm{ad}_{i_!i^!}} & & \pi_!\pi^!
(\raya) \ar@{=}[ddd] \\
\pi_! i_! (\Th (T_{\pi i}) \wedge i^*\pi^* (\raya)) \ar@{=}[d] & & \\
\pi_!i_!(i^!(\11_X)\wedge i^* (\Th (T_{\pi})\wedge \pi^*\raya)) \ar@{=}[d]  & & \\
\pi_!(i_!(i^!(\11_X))\wedge \pi^!(\raya ))\ar[rr]^-{\mathrm{ad}_{i_!i^!}}& &
\pi_!\pi^!(\raya), }
$$
where the vertical isomorphisms arise from duality and the projection formula (\cf
(\cf \cite[2.3.10]{Ayoub}).
\end{obser}

\begin{obser}\label{Obser Imagen refinada}
If $i\colon Z\to X$ is a closed immersion of codimension $d$, then the direct image
can be refined to a direct image with support. More precisely, let us write $\p_i
\colon  \EE ^{*,*}(Z)\longrightarrow \EE^{*+2d,*+d}_{Z}(X)$ for the morphism defined
as follows:
$$
\xymatrix@R=6pt{\EE (Z) \ar@{=}[d] \ar@{-->}[rr]^-{\p_i} & & \EE_Z(X) \ar@{=}[dddd] \\
\Hom_{\SH(Z)}(\11_Z , \Th(T_Z) \wedge \EE_Z) \ar@{=}[d] & & \\
\Hom_{\SH(Z)}(\11_Z , \Th(-N_{Z/X}) \wedge \Th(i^*T_X) \wedge \EE_Z) \ar[dd] & & \\
& & \\
\Hom_{\SH(Z)}(\11_Z, \Th(-N_{Z/X})\wedge \EE_Z)\ar[rr]^-{\sim} & &
\EE(\Th(N_{Z/X})),}
$$
where the right equality is a consequence of the deformation to the normal bundle
(\cf \cite[\S 1.3]{Deglise}).
\end{obser}

\begin{teo}\label{Teo Propiedades Imagen Directa}
Let $(\EE, c_1)$ be an oriented absolute ring
spectrum and $\MM$ be an absolute $\EE$-module.

The direct image in $\MM$-cohomology defined in \ref{LaDefinicion} satisfies:

\begin{enumerate}
\item
{\bf Functoriality}: for any proper morphisms $ Z \xrightarrow{g} Y \xrightarrow{f}
X$,
$$
(fg)_*=f_*g_* \ .
$$

\item
{\bf Normalization}: for any closed immersion $i\colon H\to X$ of codimension one
$$
i_*(a)=i_\flat (a \cdot c_1^H(L_H)) \ ,
$$
where $i_\flat$ is the morphism of forgetting support (\cf Remark \ref{Obser olvidar
soporte}).

\item
{\bf Key formula}: for any closed inmersion\footnote{The key formula also holds for
general projective morphisms, \cf \ref{Coro Unicidad}} $i \colon Z\to X$ of
codimension $d$
$$
\pi^* i_* (a)= j_*(c_{d-1} (K)\cdot \pi '^*(a)) \ ,
$$
where $\pi^*\colon B_ZX\to X$ is the blowing-up of $Z$ in $X$ with exceptional
divisor $j\colon \PP(N_{Z/X})\to B_ZX$.

\item
{\bf Projection formula}: for any proper morphism $f \colon Y \to X$ and any elements
$b \in \EE(X)$ and $m \in \MM(Y)$, it holds:
$$
f_*(f^* (b) \cdot m )= b \cdot f_*(m)
$$
(An analogous formula holds for $n \in \MM(X)$ and $a \in \EE(Y)$).

\end{enumerate}
\end{teo}
\demo Functoriality is evident from the definition.

Both normalization and the projection formula follow from the comparison made in
Remark \ref{Obser Gysin Thom}. Indeed, for a closed immersion $i\colon Z\to X$, the
morphism $\EE (\Th(-T_Z)) \to \EE (\Th(-T_X))$ maps $a\cdot \t^\refi (-T_Z)$ into $
i_\flat (a\cdot \t ^\refi(N_{Z/X}))\cdot \t^\refi(-T_X) $ (\cf
\cite[2.4.8]{Deglise}). Hence
\begin{equation}\label{ClaseThomYSeccionCero}
i_* (1) = i_\flat (\t ^{\refi}(N_{Z/X})) \ .
\end{equation}
In the case of codimension 1, one can explicitly compute that $\t^\refi
(N_{H/X})=c_1^H(\mathcal{I}_H^*)$, where $\mathcal{I}_H$ is the sheaf of ideals
defining $H$ (\cf \cite[2.19]{Navarro Articulo}). In addition, formula
(\ref{ClaseThomYSeccionCero}) also shows that our direct image coincides with that of
\cite[2.3.1]{Deglise}. Therefore it also satisfies the key formula (\cf
\cite[2.4.2]{Deglise}).

Finally, the projection formula is equivalent to its analogue in Borel-Moore
homology. To be precise, we have to check that for $b \in \EE^{p,q}(X)$ and $m\in
\MM^\BM_{r,s} (Y)$,
$$
f_* (f^*(b)\cdot m)=b\cdot f_*(m).
$$

Let $\pi_{\scriptscriptstyle X}\colon X\to S$ and $\pi_{\scriptscriptstyle Y}\colon
Y\to S$ be the structural morphisms. Without loss of generality, assume $p=q=r=s=0$.
By definition, the left hand side of the equation is the adjoint to the morphism
$$
\pi_{{\scriptscriptstyle X !}}\pi_{{\scriptscriptstyle X}}^* \11_S \xrightarrow{\ad}
\pi_{{\scriptscriptstyle Y !}}\pi_{{\scriptscriptstyle Y}}^* \11_S
\xrightarrow{f^*(b)} \pi_{{\scriptscriptstyle Y !}} \EE_Y\simeq \EE_S\wedge
\pi_{{\scriptscriptstyle Y !}}(\11_Y)\xrightarrow {(1\wedge m)} \EE_S\wedge \MM_S
\xrightarrow \mu \MM_S.
$$
The right hand side is the adjoint of
$$
\pi_{{\scriptscriptstyle X !}}\pi_{{\scriptscriptstyle X}}^* \11_S
\xrightarrow{\pi_{{\scriptscriptstyle X !}}(b)} \pi_{{\scriptscriptstyle X !}} \EE_X
\simeq \EE_S \wedge \pi_{{\scriptscriptstyle X !}}\pi_{{\scriptscriptstyle
X}}^*\11_S\xrightarrow{1\wedge \ad} \EE_S \wedge \pi_{{\scriptscriptstyle Y
!}}(\11_Y) \xrightarrow {1\wedge m}\EE_S\wedge \MM_S\xrightarrow \mu \MM_S.
$$

These two morphisms agree due to the same argument of \cite[1.2.10.E7]{Deglise}
replacing  inverse images by exceptional images and using the analogous
compatibility.

\qed

\begin{coro} \label{Coro Unicidad}
Let $f\colon Y\to X$ be a projective morphism and $(\EE,c_1)$ be an absolute oriented
ring spectrum.

The direct image $f_*\colon \EE(Y)\to \EE(X)$ defined in \ref{LaDefinicion} coincides
with that of \cite[3.2.6]{Deglise} and \cite[2.32]{Navarro Articulo}
\end{coro}

\demo These three definitions satisfy the conditions of the uniqueness result from
\cite[2.34]{Navarro Articulo}. \qed

\medskip
A direct image for the cohomology defined by absolute modules over absolute oriented
ring spectra was defined in \cite[\S 2]{Navarro Articulo}. For the sake of
completeness, let us prove here a uniqueness result.

\begin{teo}
Let $(\EE, c_1)$ be an oriented absolute ring spectrum and $\MM$ be an absolute
$\EE$-module.

There exists a unique way of assigning, for any projective morphism
$f\colon Y\to X$ between smooth schemes, a group morphism $f_*\colon \MM (Y) \to
\MM(X)$ satisfying the following properties:
\begin{enumerate}

\item
{\bf Functoriality}: for any projective morphisms $Z \xrightarrow{g} Y
\xrightarrow{f} X$, $(fg)_*=f_*g_*$.

\item
{\bf Normalization}: for any closed immersion $i\colon H\to X$ of codimension one,
$i_*(a)=i_\flat (a \cdot c_1^H(L_H))$.

\item
{\bf Key formula}: for any closed inmersion $i \colon Z\to X$ of
codimension $n$,
$$ \pi^* i_* (a)= j_*(c_{n-1} (K)\cdot \pi
'^*(a)) \ , $$ where $\pi^*\colon B_ZX\to X$ is the blowing-up of $Z$ in $X$ with
exceptional divisor $j\colon \PP(N_{Z/X})\to B_ZX$.

\item
{\bf Projection formula}: for any projective morphism $f \colon Y \to X$ and any
elements $b \in \EE(X)$ and $m \in \MM(Y)$, it holds:
$$
f_*(f^* (b) \cdot m )= b \cdot f_*(m) \ .
$$
\end{enumerate}
\end{teo}

\demo Any projective morphism $f\colon Y\to X$ factors as the composition of a closed
embedding $i\colon Y\to \PP^n_X$ and a projection $p\colon \PP^n_X\to X$.

Properties 1, 2 and 3 characterize the direct image for closed immersions. The
projection formula characterizes the direct image for a projection $p$ due to Lemma
\ref{Lema Fibrado Proyectivo Modulos}.

\qed

\begin{obser} \label{Obser Unic Imagen Directa}
A uniqueness result for direct images of proper morphism is also true, due to work in
progress of F. D\'{e}glise. However, such a statement requires to take into account
Borel-Moore homology and therefore the framework of \emph{bivariant theories} is more
convenient (\cf \cite{FM}).
\end{obser}

\section{Riemann-Roch theorems}

Let $(\EE , c_1)$ and $(\FF,\bar c_1)$ be oriented absolute ring spectra. To avoid
confusion, we overline notation to refer to elements and morphisms in the
$\FF$-cohomology.

\begin{parra}\label{Parra Orientaciones}
If $\varphi \colon (\EE , c_1) \to (\FF , \bar{c}_1)$ is a morphism of oriented
absolute ring spectra, then
$$
\varphi (c_1)=G(\bar c_1)\cdot \bar c_1
$$
for a unique invertible series $G(t)\in (\FF(S)[[t]])^* = \FF
(\mathbb{P}_{S}^\infty)$  (\cite[1.34]{Navarro Articulo}).

Let us write $G_\times $ for the multiplicative extension of $G(t)$. To
be more precise, if $V\to X$ is a vector bundle, the splitting principle assures the
existence of a base change $\pi\colon X'\to X$, injective in cohomology, such that
$\pi^*V=L_1+\cdots + L_r$ is a sum of line bundles in $K_0(X')$
(\cf \cite[1.35]{Navarro Articulo}); we use the notation:
$$
G_\times(V)\coloneqq G(L_1)\cdot \ldots \cdot G(L_r)\in \FF(X) \, ,
$$ where $G(L)$ stands for $G(c_1(L))$.
\end{parra}

\begin{obser}\label{Obser RR Deglise}
Let $\varphi \colon (\EE , c_1) \to (\FF , \bar{c}_1) $ be a morphism of oriented
absolute ring spectra and let $G\in \FF[S][[t]]$ be the unique series such that
$\varphi (c_{1})=G(\bar c_1)\cdot \bar c_1$.

If $i\colon Z \to X$ is a closed inmersion, and $\,N_{Z/X}\,$ stands for its normal
bundle, then D\'{e}glise (\cite[4.2.3]{Deglise}) have proved that, for any $\,a \in \EE
(Z)$, the following Riemann-Roch formula holds:
\begin{equation}\label{RRDeglise}
\varphi(\p_{i}(a))=\bar \p_i\bigl(G^{-1}_\times(-N_{Z/X})\cdot \varphi(a)\bigr) \ .
\end{equation}


\end{obser}

To be more precise, the following lemma proves the formula we will require later.

\begin{lema}
Let $\varphi\colon  (\EE , c_1) \to (\FF , \bar c_1)$ be a morphism of oriented
absolute ring spectra such that $\varphi (c_1)=G(\bar c_1)\cdot \bar c_1$.

If $ V \to X$ is a vector bundle, the following equality holds on $\FF ( \Th (V))$:
$$
\varphi (\t^\refi (V)) \, = \, G^{-1}_\times (-V)\, \bar \t^{\refi}(V)  \ .
$$
\end{lema}

\demo If $s\colon X\to V$ denotes the zero section, then $\p_s (1) = \t^\refi (V) \in
\EE (\Th V)$, (\cf equation (\ref{ClaseThomYSeccionCero}) and Remark \ref{Obser
Imagen refinada}). Moreover, $s$ is a closed inmersion whose normal bundle equals $\,
V$; thus, the Riemann-Roch formula (\ref{RRDeglise}) unfolds in this case into:
\begin{align*}
\varphi (\t^\refi (V)) &= \varphi (\p_s(1))\buildrel{(\ref{RRDeglise})}\over{=} \,
\bar \p_s \left( G^{-1}_\times (-V) \, \varphi (1) \right) = \bar \p_s \left(
G^{-1}_\times (-V) \right) = G^{-1}_\times (-V) \, \bar \t^{\refi}(V) \ .
\end{align*}

\qed

\begin{teo}\label{Teo RR motivico}
Let $f\colon Y\to  X$ be a proper morphism between smooth $S$-schemes. Let
$\varphi\colon  (\EE , c_1) \to (\FF , \bar c_1)$ be a morphism of oriented absolute
ring spectra such that $\varphi (c_1)=G(\bar c_1)\cdot \bar c_1$ and let
$G^{-1}_\times$ stand for the multiplicative extension of $G^{-1}\in \FF(S)[[t]]$.

The following square commutes:
$$
\xymatrix{ \EE(Y)\ar[r]^{f_*} \ar[d]_{G^{-1}_\times (T_Y)\varphi}& \EE(X)
\ar[d]^{G^{-1}_\times (T_X)\varphi}\\
\FF (Y) \ar[r]^{\bar f_*} & \FF (X) }
$$

\end{teo}

\demo The thesis follows from the commutativity of the following diagram: 
$$
\xymatrix@C=9pt{ \EE(Y)\ar[rr]^-{\sim}\ar[d]_{G^{-1}_\times (T_Y) \varphi} && \EE
(\Th (-T_Y))\simeq \EE^\BM (Y) \ar[d]_\varphi \ar[rr]^-{f_*} && \EE^\BM (X)= \EE (\Th
(-T_X)) \ar[d]^\varphi& & \EE(X) \ar[ll]_-{\sim} \ar[d]^-{G^{-1}_\times (T_X)
\varphi}
 \\
\FF(Y)\ar[rr]^-{\sim}& & \FF ( \Th (-T_Y))\simeq \FF ^\BM (Y) \ar[rr]^-{f_* } && \FF
^\BM(X)\simeq \FF ( \Th (-T_X)) & & \FF(X) \ar[ll]_-{\sim}. }
$$
Let us  check each of the three squares separately.

For any element $a \in \EE^{p,q}(Y)$, we have $a \cdot \t^\refi (-T_Y)\in \EE ( \Th
(-T_Y)) \simeq \EE _{-p-2n,-q-n}^\BM (Y)$. Since $\varphi$ is morphism of rings and
due to the previous Lemma,
$$
\varphi (a \cdot \t^\refi (-T_Y)) = \varphi (a) \cdot \varphi (\t^\refi(-T_Y)) =
\varphi (a) \cdot G^{-1}_\times (T_Y) \cdot \bar \t^{\refi} (-T_Y)\ .
$$
Hence, if $n=\dim_S Y$, the following diagram commutes:
$$
\xymatrix@C=12pt{ \EE^{p,q}(Y)\ar[rr]^-{\sim}\ar[d]_{G^{-1}_\times (T_Y) \varphi} &&
\EE^{p+2n,q+n} (\Th (-T_Y))\simeq \EE _{-p-2n,-q-n}^\BM (Y) \ar[d]^\varphi \\
\FF^{p,q}(Y)\ar[rr]^-{\sim}&& \FF^{p+2n,q+n} ( \Th (-T_Y)) \simeq \FF
_{-p-2n,-q-n}^\BM(Y). }
$$

\medskip

Analogously, if $d=n-\dim_S X$ the diagram
$$
\xymatrix@C=12pt{ \EE _{-p-2n,-q-n}^\BM (X) \simeq  \EE^{p+2n , q+n} ( \Th (-T_X))
\ar[d]_\varphi \ar[rr]^-\sim & & \EE^{p+2d,q+d}(X) \ar[d]^{G^{-1}_\times (T_X) \varphi } \\
\FF_{-p-2n,-q-n}^\BM(X) \simeq \FF^{p+2n , q+n} ( \Th (-T_X)) \ar[rr]^-\sim & &
\FF^{p+2d,q} (X) }
$$
commutes.

\medskip

Finally, the following square also commutes,
$$
\xymatrix@C=12pt{
\EE _{-p-2n,-q-n}^\BM (Y) \ar[d]_\varphi \ar[rrr]^-{f_*}   &&& \EE _{-p-2n,-q-n}^\BM
(X) \ar[d]^\varphi \\
\FF_{-p-2n,-q-n}^\BM(Y)\ar[rrr]^-{\bar f_*} &&&  \FF _{-p-2n,-q-n}^\BM(X) }
$$
because direct image in Borel-Moore homology is defined out of the adjunction
$f_!f^!\to \mathrm{Id}$, which is a natural transformation and hence, compatible with
morphisms of spectra.

\qed

The analogous statement for absolute modules is a consequence of the theorem above.

\begin{coro}\label{Teo RR modulos}
Let $f\colon Y\to  X$ be a proper morphism between smooth $S$-schemes.

Let $\varphi\colon (\EE , c_1) \to (\FF , \bar c_1)$ be a morphism of oriented
absolute ring spectra such that $\varphi (c_{1})=G(\bar c_1)\cdot \bar c_1\in \FF
(\PP_S ^\infty)$. Also, let $\Phi\colon \MM \to \bar{\MM}$ be a $\varphi$-morphism of
absolute modules .

The following diagram commutes:
$$
\xymatrix{ \MM(Y)\ar[r]^{f_*} \ar[d]_{G^{-1}_\times (T_Y)\Phi}& \MM(X) \ar[d]^{G^{-1}
_\times
(T_X)\Phi} \\
\bar{\MM} (Y) \ar[r]^{\bar f_*} & \bar{\MM} (X) }\phantom{G^{-1}_\times} \ .
$$

\end{coro}

\bigskip

\subsubsection*{Riemann-Roch for $K$-theory and motivic cohomology}

Let $\KGL_\QQ$ be the $\QQ$-localization ot the $K$-theory ring spectrum and $\HB$ be
the Beilinson's motivic cohomology spectrum (\cf Example \ref{Ejem espectros}).

The higher Chern character $\ch\colon \KGL_\QQ \to \bigoplus_{i}\HB [2i](i)$ is an
isomorphism of absolute ring spectra (\cite[5.3.17]{Riou3}). It does not preserve
orientations, but satisfies
$$
\ch (c_1)=\frac{1-e^{\bar{c}_1}}{\bar{c}_1} \ .
$$
Consequently, let us consider the inverse of this series, called the {\it Todd
series}:
$$
\Td (t) \coloneqq \left( \frac{1-e^t}{t}\right)^{-1} = \frac{t}{1-e^t} \ .
$$

\begin{coro}\label{Teo Riemann Roch}
Let $f\colon Y\to  X$ be a proper morphism between smooth $S$-schemes.

The following diagram, relating higher $K$-theory and motivic cohomology with
rational coefficients, commutes:
$$
\xymatrix{
K(Y)_\QQ\ar[r]^{f_*} \ar[d]_{\Td(T_Y)\ch}& K(X)_\QQ\ar[d]^{\Td(T_X)\ch}\\
H_{\M}(Y,\QQ) \ar[r]^{ f_*} & H_{\M}(X,\QQ) }
$$
\end{coro}

\begin{obser}
To our knowledge, this is the first Riemann-Roch statement involving higher
$K$-theory that does not assume any projective hypothesis (\cf \cite{Gillet},
\cite{Deglise}, \cite{HS}, \cite{Navarro Articulo}).

As regards to Grothendieck-Riemann-Roch statements---involving the $K_0$ group---the
most general formula without projective assumptions we know, due to Fulton and Gillet
(\cite{FG}),  applied to schemes defined over a field.
\end{obser}

Theorem \ref{Teo RR modulos} also applies to the modules we
described in Example \ref{Ejem espectros}. In particular, we obtain the following
result:

\begin{coro}\label{Teo RR Relativo}
Let $f\colon Y\to X$ be a proper morphism between smooth $S$-schemes. Let $g\colon T
\to X$ be a morphism of schemes and let us write $g_{\scriptscriptstyle Y}\colon
T\times_X Y\to Y$ for the induced map.

If either $g$ is proper or $f$ is smooth, then the following diagram commutes
$$
\xymatrix{ K(g_{{\scriptscriptstyle Y}})_\QQ\ar[r]^{f_*} \ar[d]_{\Td(T_Y)\ch}&
K(g)_\QQ\ar[d]^{\Td(T_X)\ch}\\
H_{\M}(g_{{\scriptscriptstyle Y}},\QQ) \ar[r]^{f_*} & H_{\M}(g,\QQ)\, . }
$$
\end{coro}

\subsubsection*{Arithmetic Riemann-Roch}

Let $\hatKGL$ and $\hatHB$ be the arithmetic $K$-theory and motivic cohomology ring
spectra defined by Holmstrom-Scholbach in \cite{HS} (see Example \ref{Ejem modulos}).
Let us also consider the arithmetic Chern character (\cf \cite[4.2]{HS}, \cite[\S
1.2]{Navarro Articulo}), which is a $\ch$-morphism of modules:
$$
\hatch\colon \hatKGL\to \hatHB \ .
$$

As a consequence of Corollary \ref{Teo RR modulos}, it also follows an arithmetic
formula: 

\begin{coro}\label{Teo RR Aritmetico}
Let $f\colon Y \to X$ be a proper morphism between smooth schemes over an arithmetic
ring.

The following diagram commutes
$$
\xymatrix{ \hatKH(Y)_\QQ\ar[r]^{f_*} \ar[d]_{\Td(T_Y)\hatch}&
\hatKH(X)_\QQ  \ar[d]^{\Td(T_X)\hatch}\\
\widehat {H}_{\M}(Y,\QQ) \ar[r]^{ f_*} & \widehat{H}_{\M}(f,\QQ) .}
$$
\end{coro}

\end{document}